\documentclass[12pt]{amsart}
\usepackage{amsmath, amssymb}

\author{}
\title{}
\date{}
\newcommand{\bbC}{{\mathbb C}}

\newcommand{\bbK}{{\mathbb K}}
\newcommand{\bbR}{{\mathbb R}}
\newcommand{\bbQ}{{\mathbb Q}}
\newcommand{\bbZ}{{\mathbb Z}}

\newcommand{\cA}{{\mathcal A}}
\newcommand{\cB}{{\mathcal B}}
\newcommand{\cC}{{\mathcal C}}
\newcommand{\cD}{{\mathcal D}}
\newcommand{\cP}{{\mathcal P}}

\newcommand{\al}{\alpha}
\newcommand{\be}{\beta}

\newcommand{\hgt}{\mbox{h}}
\newcommand{\rP}{\mbox{P}}

\newtheorem{lemma}{Lemma}
\newtheorem{theorem}{Theorem}

\begin{document}
\baselineskip=17pt

\vspace{10mm} 

\begin{center}
{\bf Primitive divisors of Lucas and Lehmer sequences, II} \\ 
\vspace{7 mm}
Paul M Voutier \\ 
Department of Mathematics \\ 
City University \\ 
Northampton Square \\ 
London, EC1V 0HB, UK 
\end{center}

\vspace{7 mm}

\begin{abstract}
Let $\al$ and $\be$ be conjugate complex algebraic integers 
which generate Lucas or Lehmer sequences. We present an algorithm 
to search for elements of such sequences which have no primitive 
divisors. We use this algorithm to prove that for all $\al$ and 
$\be$ with $\hgt(\be/\al) \leq 4$, the $n$-th element of these 
sequences has a primitive divisor for $n > 30$. In the course 
of proving this result, we give an improvement of a result of 
Stewart concerning more general sequences.  
\end{abstract}

\vspace{7 mm}



\vspace{7 mm}

\noindent
{\bf 1. Introduction}

\vspace{7 mm}

Let $\al$ and $\be$ be algebraic numbers such that 
$\al+\be$ and $\al \be$ are relatively prime non-zero 
rational integers and $\al / \be$ is not a root of unity. 
The sequence ${ \left( u_{n} \right) }_{n=0}^{\infty}$ 
defined by $\displaystyle u_{n} 
= \frac{{\al}^{n}-{\be}^{n}}{\al - \be}$ 
for $n \geq 0$ is called a {\it Lucas sequence}. 

If, instead of supposing that $\al+\be \in \bbZ$, we only 
suppose that $(\al+\be)^{2}$ is a non-zero rational integer, 
still relatively prime to $\al \be$, then we define the 
{\it Lehmer sequence} ${ \left( u_{n} \right) }_{n=0}^{\infty}$ 
associated to $\al$ and $\be$ by 
\begin{displaymath}
u_{n} = \left\{ \begin{array}{ll}
	\displaystyle \frac{{\al}^{n}-{\be}^{n}}{\al - \be}  
	& \mbox{if $n$ is odd} \\
       \displaystyle
\frac{{\al}^{n}-{\be}^{n}}{{\al}^{2}-{\be}^{2}}
	& \mbox{if $n$ is even.} \\
	       \end{array}
	\right.
\end{displaymath}

We say that a prime number $p$ is a {\it primitive divisor}
of a Lucas number $u_{n}$ if $p$ divides $u_{n}$ but does not
divide $(\al -\be)^{2}u_{2} \ldots u_{n-1}$. Similarly, $p$ is a
primitive divisor of a Lehmer number $u_{n}$ if $p$ divides $u_{n}$
but not $(\al^{2} - \be^{2})^{2}u_{3} \ldots u_{n-1}$.

Stewart \cite[p.80]{Ste1} showed, as a consequence of his Theorem 1, 
that if $n > C$ then $u_{n}$ has a primitive divisor, where 
$C=e^{452}2^{67}$ for Lucas sequences and $C=e^{452}4^{67}$ for 
Lehmer sequences. In Theorem~\ref{thm:thm2}, we shall obtain an 
improvement over Theorem 1 of \cite{Ste1} as well as decreasing 
the size of $C$. 

In an earlier article \cite{Vout1}, we enumerated all Lucas and 
Lehmer sequences whose $n$-th element has no primitive divisor 
for certain $n \leq 12$ and all $12 < n \leq 30$. We also presented 
some evidence to support the conjecture made there that for $n > 30$, 
the $n$-th element of any Lucas or Lehmer sequence always has a 
primitive divisor. 

Here, we present some further results concerning this conjecture. Our 
main result, Theorem~\ref{thm:thm1}, states that the conjecture is true 
if the absolute logarithmic height of $\al$ is small. In addition to 
providing further evidence for the validity of the conjecture (or 
at least not providing a counterexample), this result will also be 
useful in a forthcoming work where we shall make further improvements 
to the size of $C$.  

Throughout this paper, we shall denote by $\hgt(\al)$, the absolute 
logarithmic height of the algebraic number $\al$, which we shall 
define via its relationship to the minimal polynomial of $\al$ 
over $\bbZ$. Suppose that 
\begin{displaymath}
a_{d}X^{d} + a_{d-1}X^{d-1} + \ldots + a_{0} 
= a_{d} \prod_{i=1}^{d} \left( X-\al_{i} \right) \in \bbZ [X] 
\end{displaymath}
is the minimal polynomial of $\al$ over $\bbZ$ with $a_{d} > 0$, 
then we define 
\begin{displaymath}
\hgt(\al) = \frac{\log a_{d} 
		 + \sum_{i=1}^{d} \log \max \left( 1, |\al_{i}| \right) } 
		{d}.  
\end{displaymath}              

\begin{theorem}
\label{thm:thm1}
Suppose $\al$ and $\be$ generate a Lucas or Lehmer sequence with 
$\hgt(\be/\al) \leq 4$. Then, for all $n > 30$, the $n$-th element   
of this sequence has a primitive divisor. 
\end{theorem}

We prove this result by using Stewart's idea \cite[Section 5]{Ste1} 
of looking at certain Thue equations. For any Lucas or Lehmer sequence
$\left( u_{n} \right)_{n=0}^{\infty}$, there is a pair of integers $(p,q)$,
dependent only on the sequence, such that if $u_{n}$ has no primitive divisor
then $(p,q)$ is a solution of one of certain finitely many Thue 
equations associated to $n$. We use this to show that if, for  
$n > 30$, $u_{n}$ is without a primitive divisor then $n$ must 
be the denominator of a convergent in the continued-fraction 
expansion of $\arccos (p/(2q))/(2\pi)$. The advantage gained 
by this is that the convergents of real numbers grow quite quickly 
and so the problem of checking each $n$ less than $2 \cdot 10^{10}$ 
is reduced to checking no more than fifty such $n$. 

In fact, we will show that if $u_{n}$ has no primitive divisor, then  
the convergent $k/n$ must be an extremely good approximation to 
$\arccos (p/(2q))/(2\pi)$, so good that except for a few exceptional 
cases with $n$ small, we can show directly that $k/n$ is not sufficiently 
close to the number in question and therefore, eliminate $n$ from 
consideration. In the remaining cases, a direct examination of $u_{n}$ 
proves our desired result.

Stewart's upper bound for $n$, stated above, is quite large and would 
thus give rise to extremely long calculations just to determine the
convergents. Fortunately it is now possible to reduce this upper bound 
considerably. Because of its benefit to our work here, we shall determine 
such a smaller upper bound. In fact, we establish a more general result 
which is an improvement over Theorem 1 in Stewart's paper \cite{Ste1},  
whose proof requires little more effort than proving the more 
specific result which only applies to Lucas and Lehmer sequences. 

\begin{theorem}
\label{thm:thm2}
{\rm (i)} Suppose $\al$ and $\be$ are algebraic integers with 
$\be/\al$ having degree $d_{1}$ over $\bbQ$, $(\al, \be) = (1)$ 
and $\be/\al$ not a root of unity. Then there is a prime ideal 
$\cP$ which divides the ideal $(\al^{n}-\be^{n})$ but does 
not divide the ideals $(\al^{m}-\be^{m})$ for $1 \leq m < n$ 
for all 
$n > \max \left\{ 2(2^{d_{1}}-1), 4000 (d_{1} \log (3d_{1}))^{12} \right\}$. 

\noindent
{\rm (ii)} If $\al$ and $\be$ generate a Lucas or Lehmer sequence 
then the $n$-th element of this sequence has a primitive divisor 
for all $n > 2 \cdot 10^{10}$. 
\end{theorem}

\vspace{3 mm}

\noindent
{\bf 2. Preliminary Lemmas to Theorem~\ref{thm:thm2}}

\vspace{3 mm}

We shall first require a lower bound for linear forms in 
two logarithms. The work of Laurent, Mignotte and Nesterenko 
\cite{LMN} will be suitable for our needs. We also need a 
good lower bound for the height of a non-zero algebraic 
number which is not a root of unity. 

\begin{lemma}
\label{lem:lowbnd}
Suppose that $\gamma$ is a non-zero algebraic number of degree 
$D \geq 2$ over $\bbQ$ which is not a root of unity. Then 
\begin{displaymath}
\hgt(\gamma) > \frac{2}{D (\log(3D))^{3}}.  
\end{displaymath}
\end{lemma}

\begin{proof}
This is Corollary~1 of \cite{Vout2}. 
\end{proof}

\vspace{3mm}

Now let us continue. 

\begin{lemma}
\label{lem:lform}
Let $\gamma$ be a non-zero algebraic number of degree $D$ 
over $\bbQ$ which is not a root of unity and let $\log \gamma$ 
denote the principal value of its logarithm. Put 
\begin{displaymath}
\Lambda = b_{1} \log (-1) - b_{2} \log \gamma 
	= b_{1} \pi i - b_{2} \log \gamma,  
\end{displaymath}        
with $b_{1}$ a positive integer, $b_{2}$ a non-negative integer 
and $B = \max (|b_{1}|, |b_{2}|, 2)$. If $\Lambda \neq 0$ then 
\begin{displaymath}
|\Lambda| > \exp \left( -81.9 ( D \log (3D))^{3} 
			\hgt(\gamma) (\log B)^{2} \right). 
\end{displaymath}
\end{lemma}

\begin{proof}
First let us suppose that $|\gamma| \neq 1$. We can write 
$\gamma=re^{i \theta}$ where $r > 0$ and $-\pi < \theta \leq \pi$. 
Since $r^{2}=\gamma \cdot \overline{\gamma}$, we have 
\begin{displaymath}
2\hgt(r)=\hgt(r^{2})=\hgt(\gamma \cdot \overline{\gamma}) 
\leq \hgt(\gamma) + \hgt(\overline{\gamma}) = 2\hgt(\gamma), 
\end{displaymath}
so $\hgt(r) \leq \hgt(\gamma)$. Thus, by Liouville's inequality we have 
\begin{displaymath}
\left| \Lambda \right| 
  =  \left| b_{1}i \pi - b_{2}i \theta - b_{2} \log r \right| 
\geq \left| \log r \right| 
\geq 2^{-D} \exp \left( -D\hgt(\gamma) \right),  
\end{displaymath}
and the lemma follows by Lemma 1 and the fact that 
$\hgt(\gamma) \geq \log 2$ if $D=1$. 

We now turn to the case of $|\gamma|=1$. Since $\gamma$ is not a 
root of unity, $D \geq 2$. 

To obtain our lower bound for $|\Lambda|$ in this case we will use 
Th\'{e}or\`{e}me~3 of \cite{LMN}. However, this result requires that 
$b_{1}$ and $b_{2}$ be non-zero, so we must deal specially with the 
case of $b_{1}=0$.  

By Liouville's inequality and Lemma~\ref{lem:lowbnd}, 
\begin{displaymath}
|\Lambda| \geq |b_{2} \log \gamma| \geq |\log \gamma| 
	  \geq 2^{-D} \exp (-D\hgt(\gamma)) 
	  \geq \exp \left( -2D^{2} (\log (3D))^{3} \hgt(\gamma) \right).  
\end{displaymath}

It is now clear that the lemma holds in this case. 

To obtain a good constant in our lower bound we show that we may 
assume $B > 679000$. From Liouville's inequality, we obtain
\begin{displaymath}
\left| \Lambda \right| \geq 2^{-D/2} 
			    \exp \left( - \frac{DB\hgt(\gamma)}{2} \right).  
\end{displaymath}

We can use $D/2$ here instead of $D$ since $\gamma \not\in \bbR$ 
(see Exercise 3.4 of \cite{Wald}). 

So the lemma is true whenever 
\begin{displaymath}
81.9D^{2} (\log (3D))^{3} (\log B)^{2} - \frac{B}{2} 
\geq \frac{\log 2}{2\hgt(\gamma)}. 
\end{displaymath}

Since $D \geq 2$, applying Lemma~\ref{lem:lowbnd}, this inequality holds if 
\begin{displaymath}
\frac{81.9(\log B)^{2}}{B} - 0.02174 \geq \frac{0.087}{B}. 
\end{displaymath}

Using Maple, one can check that this is true for $2 \leq B \leq 679000$. 

We now invoke Th\'{e}or\`{e}me~3 of \cite{LMN}. Let 
$a=\max \left\{ 20,12.85|\log \gamma|+D\hgt(\gamma)/2 \right\}$ and 
$H=\max \left\{ 17,D\log(b_{1}/(2a)+b_{2}/(25.7\pi))/2+2.3D+3.25 \right\}$. 
Then 
\begin{equation}
\label{eq:lmn}
\log |\Lambda| \geq -9aH^{2}. 
\end{equation}

Since $a \geq 20$ and $1/(2a)+1/(25.7 \pi) < 0.0374$, 
$H \leq \max \left\{ 17, (D/2) \log B + 0.657D + 3.25 \right \}$. 
Moreover, $B > 679000$ implies that $(D/2) \log B + 0.657D + 3.25 
< 0.66994D \log B$. As this last quantity is greater than 17 
for $B > 679000$, we have $H < 0.66994D \log B$. 

We also want an upper bound for $a$ in terms of $D$ and 
$\hgt(\gamma)$. First notice that $|\log \gamma| \leq \pi$. 
Therefore, $12.85|\log \gamma|+D\hgt(\gamma)/2 
\leq D\hgt(\gamma)(40.37/(D\hgt(\gamma))+1/2)$. 
Since $D \geq 2$, we can apply Lemma~\ref{lem:lowbnd}. We obtain   
$40.37/(D\hgt(\gamma))+1/2 < 20.185 (\log (3D))^{3}+1/2 
< 20.272 (\log (3D))^{3}$. Therefore, 
$12.85|\log \gamma|+D\hgt(\gamma)/2 
< 20.272 (\log (3D))^{3} D\hgt(\gamma)$ 
for all $D \geq 2$. Moreover, this quantity is greater than 20, so 
$a < 20.272 (\log (3D))^{3} D\hgt(\gamma)$. 

Applying these estimates to (\ref{eq:lmn}), we find that our 
lemma holds. 
\end{proof}

\vspace{3 mm}        

Suppose that $\al$ and $\be$ are algebraic integers in a number field 
$\bbK$ of degree $d$ over $\bbQ$. Letting $\bbK_{1} = \bbQ (\be/\al)$, a
number field of degree $d_{1}$ over $\bbQ$, we set $\be/\al=\be_{1}/\al_{1}$,
where $\al_{1}$ and $\be_{1}$ are algebraic integers in $\bbK_{1}$ and 
$(\al_{1},\be_{1}) = \cA_{1}$. We may assume, without loss of 
generality, that $|\al_{1}| \geq |\be_{1}|$. 

We note that, unless we state otherwise, $\log z$ shall always 
denote the principal branch of the logarithmic function. 

Now let us prove: 

\begin{lemma}
\label{lem:bounds}
{\rm (i)} We have  
\begin{displaymath}
\log 2 + \log |\al_{1}| \geq \log \left| \al_{1} - \be_{1} \right| 
\geq \log |\al_{1}| - d_{1}(\hgt(\be_{1}/\al_{1}) + \log 2).  
\end{displaymath}

\noindent
{\rm (ii)} For $d_{1} \geq 2$ and $n \geq 2$, we have  
\begin{displaymath}
\log 2 + n \log |\al_{1}| 
\geq \log \left| \al_{1}^{n} - \be_{1}^{n} \right| 
\geq n \log |\al_{1}| 
     - 81.97 (d_{1} \log (3d_{1}))^{3} \hgt(\be_{1}/\al_{1}) (\log n)^{2}.  
\end{displaymath}
\end{lemma}

\begin{proof}
(i) We can write $\al_{1} - \be_{1} = \al_{1} (1 - \be_{1}/\al_{1})$. 
By Liouville's inequality, 
\begin{displaymath}
\log \left| \be_{1}/\al_{1} - 1 \right| 
\geq - d_{1}(\log 2 + \hgt(\be_{1}/\al_{1})), 
\end{displaymath}
and the result follows. 

(ii) The upper bound follows directly from the triangle inequality 
and our assumption that $|\al_{1}| \geq |\be_{1}|$. 

For the lower bound we write the quantity in question as 
\begin{displaymath}
n \log |\al_{1}| + \log \left| (\be_{1}/\al_{1})^{n} - 1 \right|.  
\end{displaymath}

Applying Lemma~2.3 of \cite{PW} with $r=1/3$ and 
$z = n \log (\be_{1}/\al_{1})$, we see that either 
\begin{displaymath}
\left| (\be_{1}/\al_{1})^{n} - 1 \right| > \frac{1}{3} \mbox{ or }
|\Lambda| = \left| n \log (\be_{1}/\al_{1}) - 2k \pi i \right| 
< 1.3 \left| (\be_{1}/\al_{1})^{n} - 1 \right| < 0.5. 
\end{displaymath}

In the first case, the lemma holds so we need only 
consider the second case. Here, we must have 
\begin{displaymath}
\left| \mbox{Im}(n \log (\be_{1}/\al_{1})) - 2k \pi \right| < 0.5. 
\end{displaymath}

Since we took the principal value of the logarithm of $\be_{1}/\al_{1}$,  
we have $-\pi < \mbox{Im}(\log (\be_{1}/\al_{1})) \leq \pi$ and so 
$|k| < n/2+0.5/(2\pi)$ or $|2k| \leq n$. 

As $\be_{1}/\al_{1}$ is, by assumption, not a root of 
unity, $\Lambda \neq 0$ and, since $n \geq 2$, we may 
apply Lemma~\ref{lem:lform} giving 
\begin{displaymath}
|\Lambda| > \exp \left( -81.9 (d_{1} \log (3d_{1}))^{3} 
			\hgt(\be_{1}/\al_{1}) 
			(\log n)^{2} \right). 
\end{displaymath}

By Lemma~\ref{lem:lowbnd}, we have 
$\log(1.3) < 0.07 ((d_{1} \log (3d_{1}))^{3} 
	     \hgt(\be_{1}/\al_{1}) \log^{2} n$,   
since $d_{1} \geq 2$ and $n \geq 2$. Our lemma follows. 
\end{proof}
			      
\begin{lemma}
\label{lem:sch}
Let $\Phi_{n}(X,Y)=Y^{\varphi(n)}\phi_{n}(X/Y)$ where $\phi_{n}(X)$ 
is the $n$-th order cyclotomic polynomial. Suppose that $\cP$ 
is a prime ideal in $\bbK$ which divides $(\Phi_{n}(\al,\be))$ for
$n > 2(2^{d_{1}}-1)$. This implies that $\cP$ divides 
$(\al^{n}-\be^{n})$. If, in addition, $\cP$ divides 
$(\al^{m}-\be^{m})$ for some $m < n$, then 
\begin{displaymath}
\mbox{\rm ord}_{\cP} \Phi_{n}(\al,\be) \leq \mbox{\rm ord}_{\cP} n. 
\end{displaymath}
\end{lemma}

\begin{proof}
This is Lemma 4 of \cite{Sch}. 
\end{proof}

\vspace{3mm}

Finally we need to bound some arithmetic functions which 
will appear throughout this article. 

\begin{lemma}
\label{lem:arith}
{\rm (i)} Let $\omega(n)$ denote the number of distinct prime factors of 
$n$. For $n \geq 3$, 
\begin{displaymath}
\omega(n) < \frac{1.3841 \log n}{\log \log n}.  
\end{displaymath}

\noindent
{\rm (ii)} For $n \geq 3$, 
\begin{displaymath}
\varphi(n) \geq \frac{n}{e^{\gamma}\log \log n + 2.50637/\log \log n}, 
\end{displaymath}
where $\gamma=0.57721\ldots$ is Euler's constant. 
\end{lemma}

\begin{proof}
(i) This follows from Th\'{e}or\`{e}me~11 of \cite{Robin}. 

(ii) This is Theorem 15 of \cite{RS}. 
\end{proof}

\vspace{3 mm}

\noindent
{\bf 3. Proof of Theorem~\ref{thm:thm2}}

\vspace{3 mm}

We may also assume that $d_{1} \geq 2$, for otherwise we can 
write $\al=\be c_{1}/c_{2}$ where $c_{1}, c_{2} \in \bbZ$ 
with $(c_{1},c_{2})=1$ and so 
$\al^{n}-\be^{n}=(c_{1}^{n}-c_{2}^{n})(\be/c_{2})^{n}$. Now 
Zsigmondy \cite{Zsig} and, independently of him, Birkhoff 
and Vandiver \cite{BV} have shown that for $n > 6$ the 
$n$-th element of such sequences always has a primitive 
divisor. 

Therefore, we may also assume that $n > 3900 (2 \log(3 \cdot 2))^{12} 
> 1.74 \cdot 10^{10}$, since Theorem~2 does not apply for smaller $n$ 
when $d_{1} \geq 2$. 

We note that 
\begin{displaymath}
\Phi_{n}(\al,\be) = \be^{\varphi(n)} \Phi_{n}(\al/\be,1) 
= \be^{\varphi(n)} \Phi_{n}(\al_{1}/\be_{1},1) 
= (\be/\be_{1})^{\varphi(n)} \Phi_{n}(\al_{1},\be_{1}).  
\end{displaymath}

Letting $\cA$ be the extension of $\cA_{1}$ in $\bbK$, we have
$\left( \be/\be_{1} \right)=\cA^{-1}$, since $(\al, \be)=(1)$, and so 
\begin{equation}
\label{eq:firsteq}
(d_{1}/d) \log \left| N_{K/Q} \left( \Phi_{n} (\al,\be) \right) \right| 
= \log \left| N_{K_{1}/Q} \left( \Phi_{n} (\al_{1},\be_{1}) \right) \right| 
  - \varphi(n) \log N_{K_{1}/Q} \left( \cA_{1} \right). 
\end{equation}

Since 
\begin{displaymath}
\Phi_{n}(\al_{1},\be_{1}) 
= \prod_{m|n} { \left( \al_{1}^{m}-\be_{1}^{m} \right) }^{\mu(n/m)}, 
\end{displaymath}
the right-hand side of (\ref{eq:firsteq}) is 
\begin{equation}
\label{eq:secondeq}
\left( \sum_{v \in M_{\infty}(K_{1})} \sum_{m|n} 
\mu(n/m) \log { \left| \al_{1}^{m} - \be_{1}^{m} \right| }_{v} \right)  
- \varphi(n) \log N_{K_{1}/Q} \left( \cA_{1} \right), 
\end{equation}
where $M_{\infty}(\bbK_{1})$ denotes the set of all archimedean 
absolute values defined on $\bbK_{1}$ up to equivalence. 

Applying Lemma~\ref{lem:bounds}, we see that the inner sum in the 
first term of this expression is at least 
\begin{eqnarray*}
&   & \log \left\{ \max \left( |\al_{1}|_{v}, |\be_{1}|_{v} \right) \right\} 
      \sum_{m|n} \mu(n/m) m 
      - \sum_{\stackrel{m|n,m>1}{\mu(n/m)=-1}} \log 2 \\ 
&   & - 81.97 (d_{1} \log (3d_{1}))^{3} \hgt(\be_{1}/\al_{1}) 
	\sum_{\stackrel{m|n,m>1}{\mu(n/m)=1}} (\log m)^{2}   
      - d_{1}(\hgt(\be_{1}/\al_{1}) + \log 2).   
\end{eqnarray*}

Combining this lower bound with 
\begin{displaymath}
\sum_{v \in M_{\infty}(K_{1})} 
\log \max \left( |\al_{1}|_{v}, |\be_{1}|_{v} \right) 
- \log N_{K_{1}/Q} \left( \cA_{1} \right) 
= \hgt(\be_{1}/\al_{1})
\end{displaymath}
and 
\begin{displaymath}
\sum_{m|n} m \mu(n/m) = \varphi(n), 
\end{displaymath}
we obtain 
\begin{eqnarray*}
(d_{1}/d) \log \left| N_{K/Q} \left( \Phi_{n} (\al,\be) \right) \right| 
& \geq & \varphi(n) \hgt(\be_{1}/\al_{1}) \\
&      & - 81.97 d_{1}^{4} (\log (3d_{1}))^{3} \hgt(\be_{1}/\al_{1}) 
	   \sum_{\stackrel{m|n}{\mu(n/m)=1}} (\log m)^{2} \\
&      & - \sum_{\stackrel{m|n}{\mu(n/m)=-1}} d_{1} \log 2 
	 - d_{1}^{2} (\hgt(\be_{1}/\al_{1}) + \log 2).  
\end{eqnarray*}

Notice that $n$ has $2^{\omega(n)-1}$ factors $m$ which satisfy 
$\mu(n/m)=1$ and the same number of factors $m$ satisfying 
$\mu(n/m)=-1$. Now, by Lemma~\ref{lem:lowbnd} and our lower 
bound for $n$,  
\begin{displaymath}
d_{1}^{2} (\hgt(\be_{1}/\al_{1}) + \log 2)  
+ \sum_{\stackrel{m|n}{\mu(n/m)=-1}} d_{1} \log 2 
< 0.005 \cdot 2^{\omega(n)} d_{1}^{4} (\log (3d_{1}))^{3} 
  \hgt(\be_{1}/\al_{1}) \log^{2} n.  
\end{displaymath}

Thus 
\begin{equation}
\label{eq:secondineq}
(d_{1}/d) \log \left| N_{K/Q} \left( \Phi_{n} (\al,\be) \right) \right| 
> \varphi(n) \hgt(\be_{1}/\al_{1}) 
  - 2^{\omega(n)} 40.99 d_{1}^{4} (\log (3d_{1}))^{3} 
  \hgt(\be_{1}/\al_{1}) (\log n)^{2},  
\end{equation}
for $d_{1} \geq 2$ and $n \geq 1.74 \cdot 10^{10}$. 

By Lemma~\ref{lem:sch}, if 
$\left| N_{K/Q} \left( \Phi_{n}(\al,\be) \right) \right| > n^{d}$,  
then there exists a prime ideal $\cP$ which divides 
$(\al^{n}-\be^{n})$ but does not divide $(\al^{m}-\be^{m})$ for any 
$m < n$. Using (\ref{eq:secondineq}) and Lemma~\ref{lem:lowbnd}, 
as well as our assumptions that $d_{1} \geq 2$ and 
$n > 1.74 \cdot 10^{10}$, this condition is satisfied if 
\begin{equation}
\label{eq:thirdineq}
\frac{\varphi(n)}{2^{\omega(n)}(\log n)^{2}} 
> 41 d_{1}^{4} (\log (3d_{1}))^{3}.  
\end{equation}

From Lemma~\ref{lem:arith}, we find that 
\begin{displaymath}
\frac{\varphi(n)}{2^{\omega(n)}(\log n)^{2}} 
> n^{0.3495}, 
\end{displaymath}
for such $n$. Therefore, (\ref{eq:thirdineq}) is satisfied for 
\begin{equation}
\label{eq:4ineq}
n > 41200 d_{1}^{11.45} (\log (3d_{1}))^{8.59}.   
\end{equation}

Since $d_{1} \geq 2$, part (i) of the theorem holds. 

(ii) Let ${ (u_{n}) }_{n=0}^{\infty}$ be a Lucas or Lehmer sequence 
generated by $\al$ and $\be$. Since $\al \be$ and $(\al + \be)^{2}$ 
are relatively prime non-zero rational integers, there exist two 
integers $p$ and $q$ such that $\al$ and $\be$ are the two roots 
of $X^{2}-\sqrt{p+2q}X+q$. Therefore, $\al, \be = 
(\sqrt{p+2q} \pm \sqrt{p-2q})/2$ and so either $\al/\be$ or $\be/\al$ 
is equal to $(p+\sqrt{p^{2}-4q^{2}})/(2q)$. Therefore we can take 
$d_{1}=2$ and so part (i) of theorem implies part (ii). 

\vspace{3 mm}

\noindent
{\bf 4. Preliminary Lemmas to Theorem~\ref{thm:thm1}}

\vspace{3 mm}

\begin{lemma}
\label{lem:cos}
Let $a$ be a non-negative real number. If $x,y \in \bbR$ with 
$-1 \leq x,y \leq 1$ and $|x-y| \leq a$ then 
\begin{displaymath}
\left| \arccos x - \arccos y \right| \leq \pi \sqrt{\frac{a}{2}}. 
\end{displaymath}
\end{lemma}

\begin{proof}
This result follows from finding the minimum value 
of the function 
\begin{displaymath}
f(x,y) = \frac{\cos x - \cos y}{(x-y)^{2}}
\end{displaymath}
on the area in $\bbR^{2}$ defined by $0 \leq x,y \leq \pi$, $x \neq y$
which is $2/\pi^{2}$ and then applying the contrapositive. 
\end{proof}

\vspace{3 mm}

Let us collect here various notations which we shall use throughout 
the remainder of this article. 

\vspace{3.0mm}

\noindent
{\bf Notations.} Given a complex-valued function $f$ defined on 
$\bbC$, we use $|f|_{1}$ to denote $\max_{|x|=1} |f(x)|$.

For a positive integer $n$, we let $g_{n}(x) \in \bbZ[x]$ be the 
minimal polynomial of $2 \cos (2 \pi/n)$ over $\bbZ$; its degree 
is $\varphi (n)/2$ if $n \geq 3$. We shall put $G_{n}(X,Y)
=Y^{\varphi(n)/2}g_{n}(X/Y)$. 

We let $m$ be the greatest odd square-free divisor of $n$. 
For such $m$, we shall write $h_{m}(X)=(X^{m}-1)/\phi_{m}(X)$.

Finally, for $n > 1$, we let $\rP(n)$ denote the largest 
prime divisor of $n$. 

\vspace{3.0mm}

As we shall see in Section 5, the crucial result needed in 
the proof of Theorem~\ref{thm:thm1} is a good lower bound 
for $|g_{n}'(2 \cos(2 \pi j/n))|$ for $(j,n)=1$. 

We will show that we need to obtain an upper bound for the 
absolute value of $h_{m}(X)$ on the unit circle which we 
find using an idea and a result of Bateman, Pomerance and 
Vaughan \cite{BPV}.  

Let us start linking these two polynomials now. 

\begin{lemma}
\label{lem:two}
Let $n \geq 3$, $1 \leq j \leq n$ with $(j,n)=1$ and 
$\zeta_{n}=\exp(2\pi i/n)$. Then 
\begin{displaymath}
\left| g_{n}'(2 \cos(2\pi j/n)) \right| 
= \left| \frac{\phi_{n}'(\zeta_{n}^{j})}{2 \sin (2\pi j/n)} \right| . 
\end{displaymath}
\end{lemma}

\begin{proof}
We can write 
\begin{eqnarray*}
\phi_{n}(X) & = & \prod_{\stackrel{1 \leq j < n/2}{(j,n)=1}} 
		  \left( X-\zeta_{n}^{j} \right) 
		  \left( X-\zeta_{n}^{-j} \right)  
	      =   \prod_{\stackrel{1 \leq j < n/2}{(j,n)=1}} 
		  \left( X^{2}-(\zeta_{n}^{j}+\zeta_{n}^{-j})X+1 \right)  \\
	    & = & \prod_{\stackrel{1 \leq j < n/2}{(j,n)=1}} 
		  \left( X^{2} + 1 - 2 \cos (2\pi j/n) X \right)  
	      =   g_{n} \left( \frac{X^{2}+1}{X} \right) X^{\varphi(n)/2}. 
\end{eqnarray*}              
		  
If $Y=(X^{2}+1)/X$ then $X=(Y \pm \sqrt{Y^{2}-4})/2=f(Y)$ and so 
\begin{eqnarray*}
g_{n}(Y)  & = & \frac{\phi_{n}(f(Y))}{f(Y)^{\varphi(n)/2}} 
		\hspace{3.0mm} \mbox{ and } \\ 
g_{n}'(Y) & = & \frac{2f(Y)\phi_{n}'(f(Y))f'(Y) 
		      -\phi_{n}(f(Y)) \varphi(n) f'(Y)} 
		     {2 f(Y)^{\varphi(n)/2+1}}.  
\end{eqnarray*}                     

Since $f(2 \cos(2\pi j/n))=\cos(2\pi j/n) \pm i \sin(2\pi j/n)$, 
\begin{displaymath}
g_{n}'(2\cos(2\pi j/n)) 
= \frac{\phi_{n}'(\cos(2\pi j/n) \pm i \sin(2\pi j/n)) f'(2 \cos(2\pi j/n))} 
       {{\left( \cos(2\pi j/n) \pm i \sin(2\pi j/n) \right) }^{\varphi(n)/2}}.  
\end{displaymath}

Notice that $f'(Y)=(1 \pm Y/\sqrt{Y^{2}-4})/2$ so that 
\begin{displaymath}
f'(2 \cos (2\pi j/n)) 
= \frac{1}{2} \left( 1 \pm i \frac{\cos (2\pi j/n)}{\sin(2\pi j/n)} \right). 
\end{displaymath}

Hence,  
\begin{displaymath}
\left| g_{n}'(2\cos(2\pi j/n)) \right|
= \frac{\left| \phi_{n}'(\zeta_{n}^{j}) \right| \sqrt{1+ \cot^{2}(2 \pi j/n)}}
       {2}
\end{displaymath}
from which the lemma follows. 
\end{proof}

\vspace{3 mm}

To work with the cyclotomic polynomials we shall need some 
relationships which they satisfy. We give these in the next lemma. 

\begin{lemma}
\label{lem:sqft}
{\rm (i)} Let $n$ be a positive integer and let $m$ be its 
greatest odd square-free divisor. We put $m'=\gcd(2,n)m$. Then  
\begin{displaymath}
\phi_{n}(X)=\phi_{m} \left( (-1)^{m'+1} X^{n/m'} \right). 
\end{displaymath}

\noindent
{\rm (ii)} Let $p$ be a prime number and $n$ any positive integer 
not divisible by $p$. Then 
\begin{displaymath}
\phi_{pn}(X) = \frac{\phi_{n}(X^{p})}{\phi_{n}(X)}. 
\end{displaymath}

\noindent
{\rm (iii)} Let $m, m'$ and $n$ be as above. We put 
$n'=n/\gcd(n,2), h_{m}(X) = (X^{m}-1)/\phi_{m}(X)$ 
and $\zeta_{n}=\exp(2 \pi i/n)$. Then, for all $j$ 
with $(j,n)=1$, we have 
\begin{displaymath}
\left| \phi_{n}'(\zeta_{n}^{j}) \right| 
= \frac{n'}{|h_{m}((-1)^{m'+1} \zeta_{m'}^{j})|}. 
\end{displaymath}
\end{lemma}

\begin{proof}
(i) This assertion follows easily from the two relations: 
\begin{displaymath}
\phi_{2t}(X) = \phi_{t}(-X) \hspace{ 5mm} \mbox{ and }  
\hspace{5 mm} \phi_{n}(X) = \phi_{m'} \left( X^{n/m'} \right),  
\end{displaymath}
which are parts (iv) and (vi) of Proposition 5.16 from Chapter 2 
of Karpilovsky's book \cite{Karp}. 

(ii) This is again from Proposition 5.16 from Chapter 2 
of \cite{Karp}. 

(iii) Applying part (i), we find that 
\begin{displaymath}
\phi_{n}'(\zeta_{n}^{j}) 
= \frac{(-1)^{m'+1} \zeta_{n}^{(n/m')-1} 
	n \phi_{m}' \left( (-1)^{m'+1} \zeta_{m'}^{j} \right)}{m'}.  
\end{displaymath}

Now $X^{m}-1=h_{m}(X) \phi_{m}(X)$ so 
$mX^{m-1}=h_{m}(X) \phi_{m}'(X)+h_{m}'(X)\phi_{m}(X)$. 
Letting $X = (-1)^{m'+1} \zeta_{m'}^{j}$, which is always a 
primitive $m$-th root of unity, we have 
$(-1)^{(m-1)(m'+1)} m \zeta_{m'}^{j(m-1)}
=h_{m}((-1)^{m'+1} \zeta_{m'}^{j}) \phi_{m}'((-1)^{m'+1} \zeta_{m'}^{j})$ 
and the result follows. 
\end{proof}

\vspace{3 mm}

We see now that we have reduced the problem of bounding $|g_{n}'|$  
from below for primitive $n$-th roots of unity to bounding $|h_{m}|_{1}$ 
from above. To deal with this new problem, we shall now use 
ideas from \cite{BPV}.  

\begin{lemma}
\label{lem:main}
Let $m=p_{1} \ldots p_{k}$ where $p_{1}, p_{2}, \ldots, p_{k}$ 
are odd primes arranged in increasing order. Then
\begin{displaymath}
{ \left| h_{m}(X) \right| }_{1} \leq 2 \prod_{i=1}^{k-1} p_{i}^{2^{k-i-1}}. 
\end{displaymath}

In fact, if $k \geq 3$ then the factor of $2$ is not needed.
\end{lemma}

\begin{proof}
From Lemma~\ref{lem:sqft}(ii),   
\begin{equation}
\label{eq:hform}
h_{m}(X) = \frac{(X^{m}-1) \phi_{p_{1} \ldots p_{k-1}}(X)}
		{\phi_{p_{1} \ldots p_{k-1}}(X^{p_{k}})} 
	 = h_{p_{1} \ldots p_{k-1}}(X^{p_{k}})  
	   \phi_{p_{1} \ldots p_{k-1}}(X). 
\end{equation}

We now use induction on $k$ to prove the lemma. 

Since $h_{1}(X)=1, h_{p_{1}}(X)=X-1$ and 
$h_{p_{1}p_{2}}(X)=(X^{p_{2}}-1) \phi_{p_{1}}(X)$, 
the lemma is true for $k \leq 2$. 

For $k=3$, we have 
${ \left| h_{p_{1}p_{2}p_{3}} \right| }_{1}  
\leq { \left| h_{p_{1}p_{2}} \right| }_{1}  
     { \left| \phi_{p_{1}p_{2}} \right| }_{1}$.  
Using the result just established for $k=2$ and a theorem 
of Carlitz \cite{Carl} which shows that 
${ \left| \phi_{p_{1}p_{2}} \right| }_{1} <p_{1}p_{2}/2$, 
we obtain 
${ \left| h_{p_{1}p_{2}p_{3}} \right| }_{1} < p_{1}^{2}p_{2}$. 
This is the desired inequality for $k=3$. 

Suppose now that the lemma holds for some $k \geq 3$. 
We apply the following estimate of Bateman, Pomerance 
and Vaughan, which follows from Theorem 1 of their paper 
\cite{BPV} and holds for $k \geq 3$, 
\begin{displaymath}
{ \left| \phi_{p_{1} \ldots p_{k}} \right| }_{1} 
< p_{k} \prod_{i=1}^{k-1} p_{i}^{2^{k-i-1}}. 
\end{displaymath}

Thus from (\ref{eq:hform}), we have 
\begin{displaymath}
{ \left| h_{p_{1} \ldots p_{k+1}} \right| }_{1}  
\leq { \left| h_{p_{1} \ldots p_{k}} \right| }_{1}  
     { \left| \phi_{p_{1} \ldots p_{k}} \right| }_{1} 
 <   \prod_{i=1}^{k-1} p_{i}^{2^{k-i-1}} 
     \times p_{k} \prod_{i=1}^{k-1} p_{i}^{2^{k-i-1}} 
 =   \prod_{i=1}^{k} p_{i}^{2^{k-i}}.  
\end{displaymath}

Hence the lemma holds. 
\end{proof}

\vspace{3.0mm}  

We need the next lemma to deal with the case $q=2$, although we will 
use it for all $q$. A simple application of the triangle inequality 
would quickly yield the inequality below with $3|q|/5$ replaced by 
$|q|/2$. However, in the case of $q=2$, this would not be sufficient 
to prove our theorem: with the lower bound that the previous lemmas 
imply for $|g_{n}'(2 \cos (2 \pi k/n))|$, the upper bound we would 
obtain for the left-hand side of (\ref{eq:satisfy}) would not decrease 
with $n$ but actually grow with $n$. To refine this trivial estimate, 
it seems we must resort to an argument like the one which follows.

\begin{lemma}
\label{lem:cosbnd}
Let $n > 30$ be a positive integer and let $p$ and $q$ be 
non-zero integers with $q \geq 2, |p| < 2q$ and 
\begin{displaymath}
\left| G_{n}(p,q) \right| \leq \rP(n/(n,3)). 
\end{displaymath}

For $1 \leq j < n/2$ with $(j,n)=1$, we put 
$\beta_{n}^{(j)}=p-2q\cos(2 \pi j/n)$. Define $k$ by 
$\displaystyle |\beta_{n}^{(k)}| 
= \min_{\stackrel{j=1 \ldots n/2}{(j,n)=1}} |\beta_{n}^{(j)}|$. Then 
\begin{displaymath}
\left| \prod_{\stackrel{j=1}{j \neq k,(j,n)=1}}^{n/2} \beta_{n}^{(j)} \right| 
> \frac{(3|q|/5)^{\varphi(n)/2} |g_{n}'(2 \cos(2\pi k/n))|}{|q|}.  
\end{displaymath}
\end{lemma}

\begin{proof}
We first divide the interval $(-2,2)$ into four subintervals and 
divide the set of integers less than $n/2$ which are relatively 
prime to $n$ into four associated subsets. Let $\cA=(-2,-1), 
\cA'= \{ m: n/3 < m < n/2, (m,n)=1 \}, \cB=(-1,0), 
\cB'= \{ m: n/4 < m < n/3, (m,n)=1 \}, \cC=(0,1), 
\cC'= \{ m: n/6 < m < n/4, (m,n)=1 \}, \cD=(1,2)$ and 
$\cD'= \{ m: 0 < m < n/6, (m,n)=1 \}$. If we let 
$\varphi(k,q,n)$ denote the number of integers in  the interval 
$(nq/k,n(q+1)/k)$ which are relatively prime to $n$ then 
$|\cA'|=\varphi(6,2,n)=\varphi(n)/2-\varphi(3,0,n), 
|\cB'|=\varphi(3,0,n)-\varphi(4,0,n), 
|\cC'|=\varphi(4,0,n)-\varphi(6,0,n)$ and 
$|\cD'|=\varphi(6,0,n)$.  

Using Theorems 5--7 of \cite{Lehmer}, we have the following inequalities 
for the cardinalities of these sets of integers: 
\begin{eqnarray}
\label{eq:lehmer}
\frac{\varphi(n) - 2^{\omega(n)}}{6}  
\leq & |\cA'| & 
\leq \frac{\varphi(n) + 2^{\omega(n)}}{6} \nonumber \\ 
\frac{\varphi(n) - 3 \cdot 2^{\omega(n)}}{12}  
\leq & |\cB'| & 
\leq \frac{\varphi(n) + 3 \cdot 2^{\omega(n)}}{12} \nonumber \\ 
\frac{\varphi(n) - 4 \cdot 2^{\omega(n)}}{12}  
\leq & |\cC'| & 
\leq \frac{\varphi(n) + 4 \cdot 2^{\omega(n)}}{12} \nonumber \\ 
\frac{\varphi(n) - 2 \cdot 2^{\omega(n)}}{6}  
\leq & |\cD'| & 
\leq \frac{\varphi(n) + 2 \cdot 2^{\omega(n)}}{6}.  
\end{eqnarray}

Let us observe that $p/q \in \cA \cup \cB \cup \cC 
\cup \cD$ and 
\begin{equation}
\label{eq:lemma10}
\left| \prod_{\stackrel{j=1}{j \neq k,(j,n)=1}}^{n/2} \beta_{n}^{(j)} \right|
= \left| \prod_{\stackrel{j=1}{j \neq k,(j,n)=1}}^{n/2} 
	 { \left( 1 - \frac{\beta_{n}^{(k)}}{\beta_{n}^{(j)}} \right) }^{-1}
  \right|
  \left| g_{n}'(2\cos(2\pi k/n)) \right|. 
\end{equation}

If $p/q \in \cA$ then $p \leq -3$, since $q \geq 2$, and either 
$k \in \cA'$ or $k$ is the largest element of $\cB'$. Thus,  
$\beta_{n}^{(j)} > 3+4\cos(2 \pi j/n)$ for each $j \in \cC' \cup 
\cD'$ and so
\begin{displaymath}
\left| \beta_{n}^{(k)} \right|
\leq c_{1} = { \left( \rP(n/(n,3)) 
	       \prod_{j \in \cC' \cup \cD'}  
		    { \left( 3+4\cos(2 \pi j/n) \right) }^{-1} 
	       \right) }^{1/(\varphi(n)/2-|\cC'|-|\cD'|)}.  
\end{displaymath} 

Combining these inequalities with (\ref{eq:lemma10}), we obtain  
\begin{displaymath}
\left| \prod_{\stackrel{j=1}{j \neq k,(j,n)=1}}^{n/2} \beta_{n}^{(j)} \right|
\geq \frac{|q|^{\varphi(n)/2-1}}{2^{|\cA'|+|\cB'|}} 
     \prod_{j \in \cC' \cup \cD'} 
	   { \left( 1 + \frac{c_{1}}{3+4\cos(2 \pi j/n)} \right) }^{-1}
     \left| g_{n}'(2\cos(2\pi k/n)) \right|. 
\end{displaymath} 

Now suppose that $p/q \in \cB$. If $\beta_{n}^{(k)} < 0$ 
then $\beta_{n}^{(j)} < 0$ for $j \in \cC' \cup \cD'$,  
so $| 1-\beta_{n}^{(k)}/\beta_{n}^{(j)}| \leq 1$ for such $j$ and 
\begin{displaymath}
\left| \prod_{\stackrel{j=1}{j \neq k,(j,n)=1}}^{n/2} \beta_{n}^{(j)} \right|
\geq \frac{|q|^{\varphi(n)/2-1}}{2^{|\cA'|+|\cB'|}} 
     \left| g_{n}'(2\cos(2\pi k/n)) \right|. 
\end{displaymath} 

Since the quantity before $|g_{n}'(2\cos(2\pi k/n))|$ on the right-hand 
side of this expression is at least as large as the similar quantity 
obtained for $p/q \in \cA$, we can ignore this case. 

If $\beta_{n}^{(k)} > 0$ then $|1-\beta_{n}^{(k)}/\beta_{n}^{(j)}| \leq 1$ 
for $j \in \cA'$ so a similar analysis to that above shows that 
\begin{displaymath}
\left| \beta_{n}^{(k)} \right|
\leq c_{2} = { \left( \rP(n/(n,3))
		      \prod_{j \in \cD'}  
		      { \left( 1 + 4\cos(2 \pi j/n) \right) }^{-1} 
	       \right) }^{1/(\varphi(n)/2-|\cD'|)} 
\end{displaymath} 
and
\begin{displaymath}
\left| \prod_{\stackrel{j=1}{j \neq k,(j,n)=1}}^{n/2} \beta_{n}^{(j)} \right|
\geq \frac{|q|^{\varphi(n)/2-1}}{2^{|\cB'|+|\cC'|}} 
     \prod_{j \in \cD'} 
	   { \left( 1 + \frac{c_{2}}{1+4\cos(2 \pi j/n)} \right) }^{-1}
  \left| g_{n}'(2\cos(2\pi k/n)) \right|. 
\end{displaymath} 

If $p/q \in \cC$ then, by the same reasoning, we obtain 
\begin{displaymath}
\left| \beta_{n}^{(k)} \right|
\leq c_{3} = { \left( \rP(n/(n,3))
		      \prod_{j \in \cA'}  
		      { \left( 1 - 4\cos(2 \pi j/n) \right) }^{-1} 
	       \right) }^{1/(\varphi(n)/2-|\cA'|)} 
\end{displaymath} 
and
\begin{displaymath}
\left| \prod_{\stackrel{j=1}{j \neq k,(j,n)=1}}^{n/2} \beta_{n}^{(j)} \right|
\geq \frac{|q|^{\varphi(n)/2-1}}{2^{|\cB'|+|\cC'|}} 
     \prod_{j \in \cA'} 
	   { \left( 1 + \frac{c_{3}}{1-4\cos(2 \pi j/n)} \right) }^{-1}
     \left| g_{n}'(2\cos(2\pi k/n)) \right|. 
\end{displaymath} 

If $p/q \in \cD$, then 
\begin{displaymath}
\left| \beta_{n}^{(k)} \right|
\leq c_{4} = { \left( \rP(n/(n,3))
		      \prod_{j \in \cA' \cup \cB'}  
		      { \left( 3 - 4\cos(2 \pi j/n) \right) }^{-1} 
	       \right) }^{1/(\varphi(n)/2-|\cA'|-|\cB'|)} 
\end{displaymath} 
and
\begin{displaymath}
\left| \prod_{\stackrel{j=1}{j \neq k,(j,n)=1}}^{n/2} \beta_{n}^{(j)} \right|
\geq \frac{|q|^{\varphi(n)/2-1}}{2^{|\cC'|+|\cD'|}} 
     \prod_{j \in \cA' \cup \cB'} 
	   { \left( 1 + \frac{c_{4}}{3-4\cos(2 \pi j/n)} \right) }^{-1}
     \left| g_{n}'(2\cos(2\pi k/n)) \right|. 
\end{displaymath} 

For $n \leq 210, n=231$ and $n=462$, we can use these estimates 
to show by direct calculation that our lemma holds. 

To deal with $n > 210$, we first show that 
$\max (c_{1},c_{2},c_{3},c_{4}) < 1$ for such $n$. Using the 
above expressions for these quantities we see that this holds 
if $\min (3^{|\cA'|}, 3^{|\cD'|}) > n$. By 
(\ref{eq:lehmer}), both $|\cA'|$ and $|\cD'|$ are at 
least $\varphi(n)/6-2^{\omega(n)}/3$, so we need only prove that 
$(\varphi(n) - 2 \cdot 2^{\omega(n)}) (\log 3) > 6\log n$ 
for $n > 210$.

For $210 < n < 330=2 \cdot 3 \cdot 5 \cdot 11$, $2^{\omega(n)} \leq 8 
< n^{0.389}$. Lemma~\ref{lem:arith}(ii) yields the lower bound 
$\varphi(n) > n^{0.719}$ for $n \geq 210$. Since $\log n < n^{0.314}$ 
for $n > 210$, we need only show that $n^{0.075}(n^{0.33}-2) \log 3 > 6$ 
for $n$ in this range. But this is easily seen to be true. 

For $330 \leq n < 2310=2 \cdot 3 \cdot 5 \cdot 7 \cdot 11$, 
$2^{\omega(n)} \leq 16 < n^{0.48}$. Moreover, by Lemma~\ref{lem:arith}(i), 
for $n \geq 2310$, $2^{\omega(n)} < n^{0.9594/\log \log n} < n^{0.47}$.  
Therefore, for $n \geq 330$, $2^{\omega(n)} < n^{0.48}$. Applying 
Lemma~\ref{lem:arith}(ii) again, we find that $\varphi(n) > n^{0.73}$ 
for $n \geq 330$. Since $\log n < n^{0.31}$ for $n \geq 330$, we 
need only show that $n^{0.17}(n^{0.25}-2) \log 3 > 6$ for $n$ in 
this range which is also easily seen to be true. Therefore, 
$\max (c_{1},c_{2},c_{3},c_{4}) < 1$.

So, from our lower bounds for the absolute values of the products 
of the $\be_{n}^{(j)}$'s given above, to prove the lemma we need 
to show that 
\begin{displaymath}
\max \left( 2^{|\cA'|+|\cB'|} (4/3)^{|\cC'|+|\cD'|}, 
2^{|\cB'|+|\cC'|} (4/3)^{|\cD'|}, 
2^{|\cB'|+|\cC'|} (4/3)^{|\cA'|}, 
2^{|\cC'|+|\cD'|} (4/3)^{|\cA'|+|\cB'|} \right) 
\end{displaymath}
is less than $(5/3)^{\varphi(n)/2}$. 

Let us first show that $2^{|\cA'|}(2/3)^{|\cC'|} \geq 1$ 
and $2^{|\cD'|}(2/3)^{|\cB'|} \geq 1$. These inequalities  
will show that either the first or the last terms give the maximum 
in this expression. 

For the first of these two inequalities to be true, by (\ref{eq:lehmer}) 
we need to show that $0.08 \varphi(n) - 0.26 \cdot 2^{\omega(n)} \geq 0$. 
Similarly, the second inequality requires that the stronger inequality 
$0.08 \varphi(n) - 0.34 \cdot 2^{\omega(n)} \geq 0$ holds. So we need 
only consider this last inequality which we shall rewrite in the form 
$0.08/0.34 \geq 2^{\omega(n)}/\varphi(n)$. 

For $210 < n < 330$, we saw in a previous paragraph that 
$2^{\omega(n)}/\varphi(n) < n^{-0.33} < 0.171 < 0.08/0.34$. We also 
saw that $2^{\omega(n)}/\varphi(n) < n^{-0.25} < 0.235 < 0.08/0.34$ 
for $n \geq 330$. Therefore, our desired inequalities holds and we 
need only try to bound 
$2^{|\cA'|+|\cB'|} (4/3)^{|\cC'|+|\cD'|}$ and 
$2^{|\cC'|+|\cD'|} (4/3)^{|\cA'|+|\cB'|}$ from above.  

Notice that $|\cA' \cup \cB'|=\varphi(4,0,n)$ and that 
$|\cC' \cup \cD'| = \varphi(4,1,n)$. Lehmer \cite[p. 351]{Lehmer} 
has noted that $E(4,1,n)=-E(4,0,n)$, where $E(k,q,n)$ denotes 
$\varphi(n)-k\varphi(k,q,n)$, so we need only examine
\begin{displaymath}
(8/3)^{\varphi(n)/4} (3/2)^{|E(4,0,n)|/4}.
\end{displaymath}

Lehmer also gives precise information about $E(4,0,n)$ in Theorem 
6 of \cite{Lehmer}. If $n > 4$ and $4$ divides $n$ or $n$ is divisible 
by a prime congruent to 1 mod 4 then $E(4,0,n)=0$ and our proof 
is complete. If neither of these conditions is true then 
$|E(4,0,n)|=2^{\omega(n')}$ where $n'$ is as in the statement 
of Lemma~\ref{lem:sqft}(iii). Notice that when $E(4,0,n) \neq 0$, 
$n$ is not congruent to 0 mod 4, so $n'$ is the odd part of $n$. 

A direct calculation shows that for $210 < n < 750$, with 
the exceptions of $n=231$ and 462 which we considered 
above, $|E(4,0,n)|/\varphi(n) < 0.05$. Therefore, 
$(8/3)^{\varphi(n)/4} (3/2)^{|E(4,0,n)|/4} < (5/3)^{\varphi(n)/2}$ 
for $210 < n < 750, n \neq 231, 462$. Recalling that we showed by 
calculation that the lemma holds for $30 < n \leq 210$, for $n=231$ 
and for $n=462$, we now know that the lemma holds for $30 < n < 750$.

Notice that if $n < 4389 = 3 \cdot 7 \cdot 11 \cdot 19$ then either 
$|E(4,0,n)|=0$ or $2^{\omega(n')} \leq 8$, since in the latter 
case $n'$ is odd and without prime divisors congruent to 1 mod 4. 
Using Lemma~\ref{lem:arith}(ii), $\varphi(n) \geq 160$ and so 
$|E(4,0,n)|/\varphi(n) < 0.05$ for $n \geq 750$ and our lemma 
holds for $30 < n < 4389$. 

Applying the inequality $2^{\omega(n)} < n^{0.96/\log \log n}$, 
which follows from Lemma~\ref{lem:arith}(i), and part (ii) of this  
same lemma, we find that 
\begin{displaymath}
\frac{2^{\omega(n)}/4}{\varphi(n)/4} 
< \frac{n^{0.96/\log \log n} (1.7811 \log \log n + 2.51/\log \log n)}{n}. 
\end{displaymath}

The right-hand side is a monotone-decreasing function for $n \geq 10$ 
and so it is less than $0.05$ for $n \geq 4389$. Therefore
$(8/3)^{\varphi(n)/4} (3/2)^{|E(4,0,n)|/4} < (5/3)^{\varphi(n)/2}$ for 
$n \geq 4389$, which shows that the lemma is true.
\end{proof}

\vspace{3 mm}

\noindent
{\bf 5. Proof of Theorem~\ref{thm:thm1}}

\vspace{3 mm}

Let ${ (u_{n}) }_{n=0}^{\infty}$ be a Lucas or Lehmer sequence 
generated by $\al$ and $\be$. As noted in the proof of 
Theorem~\ref{thm:thm2}, there exist two integers $p$ and $q$ such 
that $\al$ and $\be$ are the two roots of $X^{2}-\sqrt{p+2q}X+q$. 
Notice that the $n$-th element of the sequence generated by $i \al$ 
and $i \be$ is just $\pm u_{n}$. Therefore, we can assume that 
$q=\al \be$ is positive. Also notice $|p| < 2q$ for otherwise 
$\al$ and $\be$ are real and Carmichael \cite{Carm}, Ward \cite{Ward} 
and Durst \cite{Durst} have shown that in this case the $n$-th element 
of these sequences has a primitive divisor for $n > 12$. 

Let us define the $\beta_{n}^{(j)}$'s and $\beta_{n}^{(k)}$ as in  
Lemma~\ref{lem:cosbnd}. Stewart \cite[Section 5]{Ste1} has shown that 
if the $n$-th element of this sequence has no primitive divisor then 
\begin{equation}
\label{eq:stewart}
\left| G_{n}(p,q) \right| 
= \prod_{\stackrel{1 \leq j \leq n/2}{(j,n)=1}} \left| \beta_{n}^{(j)} \right| 
\leq \rP(n/(n,3)) \mbox{ for $n > 12$}.
\end{equation}

Since $|p|<2q$, upon applying Lemma~\ref{lem:cosbnd}, we obtain 
\begin{displaymath}
\left| \beta_{n}^{(k)} \right| 
\leq \frac{\rP(n/(3,n))}
     {\displaystyle \prod_{\stackrel{1 \leq j \leq n/2}{j \neq k, (j,n)=1}} 
\left| \beta_{n}^{(j)} \right|} 
  <  \frac{(5/3)^{\varphi(n)/2} \rP(n/(n,3))}
	  {|g_{n}'(2 \cos(2 \pi k/n))| |q|^{\varphi(n)/2-1}},
\end{displaymath}
for $n > 30$.

Therefore, if we can show that 
\begin{equation}
\label{eq:satisfy}
\left| \frac{p}{q} - 2 \cos \left( \frac{2 \pi k}{n} \right) \right| 
< \frac{(5/(3|q|))^{\varphi(n)/2} \rP(n/(n,3))}{|g_{n}'(2 \cos(2 \pi k/n))|} 
< \frac{4}{n^{4}}, 
\end{equation}
then, by Lemma~\ref{lem:cos},  
\begin{displaymath}
\left| \frac{1}{2 \pi} \arccos \left( \frac{p}{2q} \right) 
       - \frac{k}{n} \right| 
< \frac{1}{2n^{2}}, 
\end{displaymath}
and so, by Theorem~184 of \cite{HW}, $k/n$ must be a convergent 
in the continued-fraction expansion of $\arccos (p/(2q))/(2\pi)$.

Hence we first want to show that for $n$ sufficiently large, 
the right-hand inequality of (\ref{eq:satisfy}) holds. We 
start by considering the case of $q=2$, as this is the most 
difficult one. 

\vspace{3 mm}

\noindent
{\bf 5.1. The case $q=2$}

\vspace{3 mm}

Using the notation of Lemmas~\ref{lem:sqft} and \ref{lem:main}, 
we find, from Lemma~\ref{lem:main}, that  
\begin{displaymath}
|h_{m}(X)|_{1} \leq m^{2^{k-1}/k} \leq n^{2^{\omega(n)-1}/\omega(n)},   
\end{displaymath}
for $m > 1$. 

If $m=1$, but $n > 1$, we have 
$|h_{m}(X)|_{1} = 1 \leq n^{2^{\omega(n)-1}/\omega(n)}$. 

Combining this upper bound with Lemmas \ref{lem:two} and 
\ref{lem:sqft}(iii), we obtain  
\begin{displaymath}
|g_{n}'(2 \cos (2 \pi k/n))| > \frac{n}{4 n^{2^{\omega(n)-1}/\omega(n)}},   
\end{displaymath}
for $n > 1$.

Applying this lower bound to the right-hand inequality of 
(\ref{eq:satisfy}) and squaring both sides, we want to show that 
\begin{displaymath}
(5/6)^{\varphi(n)} n^{2^{\omega(n)}/\omega(n)} \leq \frac{1}{n^{8}},  
\end{displaymath}
for $n > 30$.

To prove that this holds for $n$ sufficiently large, we take 
the logarithm of both sides, which yields 
\begin{displaymath}
\varphi(n) \log (5/6) + 2^{\omega(n)} (\log n)/ \omega(n) + 8 \log n \leq 0. 
\end{displaymath}

From Lemma~\ref{lem:arith}(ii), we see that $\varphi(n) > n^{0.8043}$ 
for $n \geq 3500$, while for the term involving $\omega(n)$ we use the 
fact that $2^{x}/x$ is a monotone increasing function for $x > 1/\log 2$,  
$2^{1}/1 = 2^{2}/2$ and Lemma~\ref{lem:arith}(i). In this manner, our 
problem is to show that 
\begin{displaymath}
-0.182n^{0.8043} + \frac{n^{0.96/\log \log n} \log \log n}{1.384} 
+ 8 \log n \leq 0. 
\end{displaymath}

For $n \geq 3500$, the sum of the second and third terms is at most 
$n^{0.5952}$. Therefore, we need only show that $-0.182n^{0.209}+1 \leq 0$, 
but this is easily seen to be true for $n \geq 3500$. So we have an initial 
bound of intermediate size. 

Notice though that we did not make full use of the Lemma~\ref{lem:main} 
in this argument. A direct calculation on a computer using the result 
given in Lemmas~\ref{lem:two},\ref{lem:sqft}(iii) and \ref{lem:main} 
shows that the right-hand inequality of (\ref{eq:satisfy}) holds for 
all $n > 1260$ when $q=2$. 

In the case of $q=2$, Lucas and Lehmer sequences can result from 
$p=-3,-1,1$ and 3. Since $G_{n}(p,q)$ is a product of terms of the 
form $p-2q\cos(2\pi i/n)$, it is quite easy to calculate $G_{n}(p,q)$, 
although care must be taken to maintain sufficient accuracy, and so 
we can check whether $u_{n}$ has primitive divisors by means of 
(\ref{eq:stewart}). However, to check $u_{n}$ for each $n$ up to 
1260 in this manner is quite time-consuming. Fortunately, one can 
quickly extract still more information from (\ref{eq:satisfy}). 
Given $n,p$ and $q$, it is easy to find the integer $k$ with 
$(k,n)=1$ which minimizes the far left-hand side of (\ref{eq:satisfy}). 
As when considering $1260 < n < 3500$, we can bound from above the 
middle quantity in (\ref{eq:satisfy}). For $q=2,p=-3,-1,1,3$ and 
$330 < n \leq 1260$, we can verify in this way that the left-hand 
inequality in (\ref{eq:satisfy}) is violated and so for such $n$, 
the $n$-th element of these sequences has a primitive divisor. But 
we still need to consider $30 < n \leq 330$. For these $n$, we use 
(\ref{eq:stewart}) as described earlier in this paragraph. 

For $n > 1260$, we have seen that $n$ must be the denominator 
of a convergent in the continued-fraction expansion for 
$\arccos(p/(2q))/(2\pi)$. 

The question arises of how to deal with these $n$. We are fortunate 
that in these cases the middle quantity in (\ref{eq:satisfy}) is 
extremely small. For such $n$, we proceed in the same manner that 
we checked the left-hand inequality in (\ref{eq:satisfy}) holds 
for $330 < n \leq 1260$, except that now we know $k$ too. 
Theorem~\ref{thm:thm2} tells us that we need only check those 
convergents $k/n$ with $n \leq 2 \cdot 10^{10}$. For each convergent 
computed with $n \leq 2 \cdot 10^{10}$, $|p/q-2\cos(2\pi k/n)|$ was 
considerably larger than the bound that the left-hand inequality 
of (\ref{eq:satisfy}) requires if $u_{n}$ were to be without a 
primitive divisor. In Table~1, for $p=-3$, we list the convergents 
with $n > 1260$ and give the logarithms of the required and actual 
bounds, denoted $d_{\mbox{req}}$ and $d_{\mbox{act}}$, respectively. 
The value of $\log |d_{\mbox{req}}|$ given in Table~1 is truncated 
to its integer part, whereas the value of $\log |d_{\mbox{act}}|$ 
is truncated to one decimal place.

Proceeding in this same way for $p=-1,1$ and 3, we are able to 
conclude that if $(u_{n})_{n=0}^{\infty}$ is the Lucas or Lehmer 
sequence generated by any of the pairs $(\al,\be) = (1 \pm \sqrt{-7})/2, 
(\sqrt{3} \pm \sqrt{-5})/2, (\sqrt{5} \pm \sqrt{-3})/2$ or 
$(\sqrt{7} \pm \sqrt{-1})/2$, then $u_{n}$ has a primitive divisor 
for $n > 30$. 

\begin{table}
\centering
\begin{tabular}{||r|r|r|r||}                                                               \hline
     $k$        &         $n$       & $\log |d_{\mbox{req}}|$ & $\log |d_{\mbox{act}}|$ \\ \hline
	    497 &              1291 &                   -116. &  -12.6                  \\ \hline
	    579 &              1504 &                    -68. &  -13.7                  \\ \hline
	   1655 &              4299 &                   -260. &  -15.4                  \\ \hline
	   3889 &             10102 &                   -459. &  -18.9                  \\ \hline
	  52212 &            135625 &                  -8207. &  -22.1                  \\ \hline
	  56101 &            145727 &                 -12970. &  -22.4                  \\ \hline
	 108313 &            281352 &                  -8086. &  -24.3                  \\ \hline
	 381040 &            989783 &                 -90228. &  -26.1                  \\ \hline
	 489353 &           1271135 &                 -90181. &  -26.7                  \\ \hline
	 870393 &           2260918 &                 -93683. &  -28.3                  \\ \hline
	2230139 &           5792971 &                -493472. &  -29.5                  \\ \hline
	3100532 &           8053889 &                -734197. &  -30.8                  \\ \hline
	8431203 &          21900749 &               -1745895. &  -32.3                  \\ \hline
       11531735 &          29954638 &               -1165244. &  -33.1                  \\ \hline
       19962938 &          51855387 &               -3104401. &  -34.1                  \\ \hline
       31494673 &          81810025 &               -5404005. &  -35.2                  \\ \hline
       51457611 &         133665412 &               -5943915. &  -35.8                  \\ \hline
       82952284 &         215475437 &              -19412834. &  -38.5                  \\ \hline
      798028167 &        2072944345 &             -144472147. &  -41.8                  \\ \hline
     1679008618 &        4361364127 &             -374075698. &  -42.8                  \\ \hline
     2477036785 &        6434308472 &             -293278284. &  -44.3                  \\ \hline
     6633082188 &       17229981071 &            -1438733756. &  -45.4                  \\ \hline
\end{tabular}
\caption{$(p,q)=(-3,2)$ Verification}
\end{table}

\vspace{3 mm}

\noindent
{\bf 5.2. The case of $q>2$}

\vspace{3 mm}

For such pairs $(p,q)$, we proceed along the same lines. The 
only difference is that less work is required for small $n$. 
We already know that the right-hand inequality of (\ref{eq:satisfy}) 
is satisfied for $n > 1260$ by our work in the previous section. 
We can check directly, as with $1260 < n < 3500$ for $q=2$, 
that the right-hand inequality of (\ref{eq:satisfy}) holds 
for $n_{q-1} \geq n > n_{q}$ where $n_{q}$ is given in Table~2.

As in the case of $30 < n \leq 330$ for $q=2$, we directly check 
those $u_{n}$ with $30 < n \leq n_{q}$ for primitive divisors and 
for larger $n$ we compare the required and actual differences of 
$|p/q-2 \cos(2\pi k/n)|$ in (\ref{eq:satisfy}) to establish our 
result. The actual difference is less than the required difference 
for all $n_{q} < n \leq 2 \cdot 10^{10}$ and all $3 \leq q \leq 3000$ 
(this corresponds to all pairs of $\al$ and $\be$ with $\hgt(\be/\al)
\leq 4$). By Theorem~\ref{thm:thm2}, Theorem~\ref{thm:thm1} now follows. 

All the calculations in this article were performed using Release 3 
of Maple V and UBASIC 8.74 on an IBM-compatible PC with an 486DX2 
running at 66 MHz. In total, the calculations required just over 100 
hours on this machine. Many of the calculations were performed using  
both systems to provide a check on the quantities obtained and the 
results were always identical up to the specified accuracy. 

\begin{table}
\centering
\begin{tabular}{||r|r||}    \hline
$q$           & $n_{q}$  \\ \hline
2             &  1260    \\ \hline
3             &   330    \\ \hline
4             &   210    \\ \hline
5             &   120    \\ \hline
6             &    90    \\ \hline
7             &    78    \\ \hline
8             &    66    \\ \hline
9,10,11       &    60    \\ \hline
12,\ldots,20  &    42    \\ \hline
$\geq 21$     &    30    \\ \hline
\end{tabular}
\caption{Values of $n_{q}$}
\end{table}

\setlength{\baselineskip}{5 mm}

\end{document}